\documentclass[10pt]{amsart}
\usepackage[cp1251]{inputenc}
\usepackage[english,russian]{babel}
\usepackage{amsmath}
\usepackage{amssymb}
\usepackage{amsfonts}

\usepackage{graphicx}
\newtheorem{lemma}{Lemma}
\newtheorem{theorem}{Theorem}

\newtheorem{remark}{Remark}

\def\Z{\mathbb{Z}}    
\def\P{{\mathbf P}}
\def\E{{\mathbf E}}

\def\logo{{\bf\huge S\raisebox{0.2ex}{\hspace{0.55ex}\raisebox{0.05ex}e\hspace{-1.65ex}$\bigcirc$}MR}}

\def\semrtop
     {
  \vbox{
     \noindent\logo
     \hspace{80mm}\raisebox{1ex}{ISSN 1813-3304 }

     \vspace{5mm}
\begin{center}
     {\large Siberian Electronic Mathematical Reports} \\[1mm]
     {\LARGE\tt{http://semr.math.nsc.ru}}\\[0.5mm]
     \end{center}
     \vspace{-3mm}
     \noindent
     \begin{tabular}{c}
     \hphantom{aaaaaaaaaaaaaaaaaaaaaaaaaaaaaaaaaaaaaaaaaaaaaaaaaaaaaaaaaaaaaaaaaaaaaa} \\
     \hline\hline
     \end{tabular}

  }
}

\begin{document}

\renewcommand{\refname}{References}
\renewcommand{\proofname}{Proof.}
\thispagestyle{empty}

\title[Random Multiple Access with energy harvesting]{Stability and instability of a random multiple access
model with adaptive energy harvesting}
\author{{S. FOSS, D. KIM, A. TURLIKOV}}%
\address{Sergey Foss
\newline\hphantom{iii} S.L. Sobolev Institute of Mathematics and Novosibirsk State University,
\newline\hphantom{iii} 630090, Novosibirsk, Russia;
\newline\hphantom{iii} Heriot-Watt University,
\newline\hphantom{iii} EH14 4AS, Edinburgh, UK}
\email{foss@math.nsc.ru, s.foss@hw.ac.uk}%

\address{Dmitriy Kim
\newline\hphantom{iii} Kazakh National Research Technical University named after K.~I.~Satpaev
\newline\hphantom{iii} and  ``EcoRisk'' LLP,
\newline\hphantom{iii} Satpaev str. 22,
\newline\hphantom{iii} 050013, Almaty, Kazakhstan}%
\email{kdk26@mail.ru}%

\address{Andrey Turlikov
\newline\hphantom{iii} St.-Petersburg University of Aerospace Instrumentation,
\newline\hphantom{iii} B. Morskaya 67,
\newline\hphantom{iii} 190000, St.-Petersburg, Russia,}%
\email{turlikov@vu.spb.ru}%

\thanks{\sc Foss, S., Kim, D., Turlikov, A.,
Stability and instability of a random multiple access model with
adaptive energy harvesting}
\thanks{\copyright \ 2016 Foss, S., Kim, D., Turlikov, A}
\thanks{\rm The authors thank the  grant 0770/GF3 of Kazakhstan
Ministry of Education and Science. }

\semrtop \vspace{1cm} \maketitle {\small
\begin{quote}
\noindent{\sc Abstract. } We introduce a model for  the classical
synchronised multiple access system with a single transmission
channel and a randomised transmission protocol (ALOHA). We assume
in addition that there is an energy harvesting mechanism, and any
message transmission requires a unit of energy. Units of energy
arrive  randomly and independently of anything else. We analyse
stability and instability conditions for this model.

 \medskip

\noindent{\bf Keywords:} random multiple access;  stochastic
energy harvesting; (in)stability; ALOHA algorithm;  generalised
Foster criterion.
 \end{quote}
}

\section{Introduction}

Nowadays the idea of usage of  energy harvesting in communication
systems is of a common interest, due to many applications. For
example, sensor networks with rechargeable batteries that are
harvesting energy from the environment can significantly extend
the lifetime of the system. Another example of using energy
harvesting was introduced in \cite{XRC}, where the authors
considered a multi-user system with the base station that employs
technology OFDM and a wireless channel to transmit both data and
energy to its users.
 Usage of an energy harvesting  mechanism in
systems with the random multiple access  presents new challenges
and, particularly, in determining their stability regions.  In
\cite{Ephr}, the authors considered a model with a decentralised
energy harvesting mechanism, with assuming individual power
supplies for a finite number of transmitting nodes, and studied
their stability properties. In this paper, we consider another
model (see fig. 1) with infinite number of transmitters and
separate power supplies per each transmitting node (with certain
limitations on the intensity of supply) and examine its stability
and instability. Earlier (see \cite{sib} and fig. 2), we
considered random multiple access models with a common
(centralised) renewable power supply and infinite number of
transmitters. We proved that an additional energy limitation on a
common energy supply may stabilise the system.  At the same time,
the model from \cite{sib} is simplistic and far from the practice
(in comparison with that from \cite{Ephr}). On the other hand, if
one considers infinite number of users in the model of \cite{Ephr},
the system may lose stability.
\begin{otherlanguage}{english}

\begin{figure}
\label{ut}
 \includegraphics[width=0.4\textwidth]{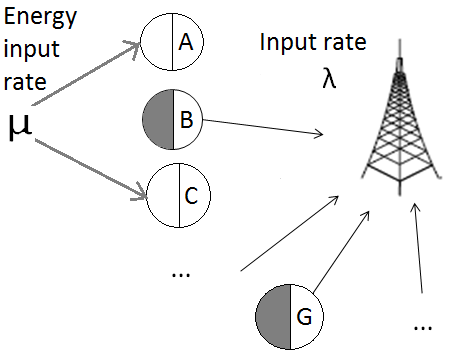}
\caption{Random access of infinite number of users, where each of
them is equipped with a battery for storing energy.}
\end{figure}

\begin{figure}
\label{jn}
\includegraphics[width=0.4\textwidth]{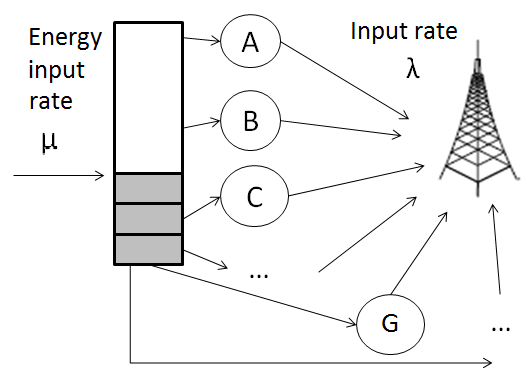}
\caption{Random access of infinite number of users with common
energy supply.}
\end{figure}

\end{otherlanguage}

 Recall that
 the classical ALOHA (see, e.g., \cite{Abr}) algorithm
 is always unstable in the system with
infinitely many users, and has a non-empty stability region
 if the number
of users is finite.

In this paper, we study a model with infinite number of
transmitters/users having individual power supplies, and with
energy harvesting mechanism with intensity that depends on the
number of active transmitters. We are interested in conditions on the
energy harvesting algorithm for (in)stability of the model.

\section{The model}
This is a system with a single transmission channel and messages that arrive and then
depart after a successful transmission.
Time is slotted, with $\xi_n$ being the number of messages that
arrive within time slot $[n-1,n)$ (the $n$-th time slot). We
assume $\{\xi_n\}$ to form an i.i.d. sequence  with a general distribution on
non-negative integers with finite mean
$\lambda\in(0,1)$. Every message is equipped with a battery of a unit volume
for storing energy.  Each 
message arrives into the system with empty battery and  waits for
energy harvesting (another model where messages arrive
into the system with units of energy was analysed in \cite{FKT}). The
harvesting mechanism within any time slot is as follows: if a
message does not have energy for transmission, its battery receives
a unit of energy with probability $\mu>0$ or nothing
otherwise,
independently of everything else. If the
battery is full, no new units of energy are accepted further
(equivalently, an energy unit may arrive, but then is rejected
by the user since the battery is already full).

 Within each time slot, any message with energy
is transmitted
 with probability $p$
(where $p\in (0,1]$ is fixed) and stays silent with probability
$1-p$. If there is a transmission of only one message, it is successful, and the message leaves
the system. If two or more messages are transmitted simultaneously, a collision occurs,
all messages stay in the system, but loose their units of energy (in other words, each unit
of energy is used for a single transmission attempt only).

The dynamics of the model  may be described as follows. Assume
that at the  start of, say, the $(n+1)$st time slot, there are
$\xi_n$ ``new'' messages that just arrived and $q_n$ ``old''
messages that arrived earlier. Let $v_n\le q_n$ be the number of
old messages that have energy units for transmission. Let
$\{U_{n,i}, -\infty <n<\infty, i\ge 1\}$ and $\{\widehat{U}_{n,i},
-\infty <n<\infty, i\ge 1\}$ be two independent families of i.i.d.
random variables that are uniformly distributed in the interval
$(0,1)$. Let ${\mathbf I} (A)$ be the indicator function of event
$A$: it takes value 1 if the event occurs, and value 0 otherwise.
Then $B_n(k,p) = \sum_{i=1}^k {\mathbf I} (U_{n,i}<p),
-\infty<n<\infty, k\ge 0$ and $\widehat{B}_n(k,\mu ) =
\sum_{i=1}^k {\mathbf I} (\widehat{U}_{n,i}< \mu ), -\infty
<n<\infty, k\ge0$ are two mutually independent families of random
variables with binomial distributions that do not depend on
anything else. We will also use random variables $D_n(k,1-p):= k -
B_n(k,p) = \sum_{i=1}^k {\mathbf I} (U_{n,i}>p)$ that have a
binomial distribution with parameters $k$ and $1-p$.

The pairs $\left(q_n,v_n\right)$ form a two-dimensional
time-homogeneous Markov chain with the following dynamics:
\begin{eqnarray*}
q_{n+1} &=& q_n - {\bf I} (B_n(v_n,p)=1) + \xi_{n},\\
v_{n+1} &=& v_n - B_n(v_n,p) + \widehat{B}_n(q_n - v_n+\xi_n,\mu).
\end{eqnarray*}
Note that we can replace $v_n-B_n(v_n,p)$ by $D_n(v_n,p)$ in the
last line.

In practice, the intensity $\mu$ of energy
harvesting is {\it adaptive}: it depends on the number $q$ of the messages in the system,
it decreases when $q$ increases, see
e.g. \cite{Can}, \cite{Tong}. In this paper, we assume that the
dependence is inverse proportional: for $n=0,1,\ldots$, the intensity $\mu = \mu_n$
in the $n$th time slot is given by
\begin{equation}\label{star}
\mu_n = \mu(q_n) = \frac{c}{q_n},
\end{equation}
for some $c>0.$

\section{Main results}

We say that the system is {\it stable} if the underlying Markov
chain is positive recurrent, and {\it unstable} if
$q_n+v_n\to\infty$ in probability, as $n\to\infty$.

To analyse the (in)stability of $(q_n,v_n)$, we introduce an
auxiliary one-dimensional Markov chain $\widetilde{V}_n$ by the following
recursion:
\begin{eqnarray}
\label{aux} \widetilde{V}_{n+1} = \widetilde{V}_n -
B_n(\widetilde{V}_n,p) + \eta_n,
\end{eqnarray}
where $\{\eta_n\}_{n=0}^\infty$ is an i.i.d. sequence of Poisson
random variables with parameter $c$. One can easily apply the
classical Foster criterion and condition ${\mathbf P} (\eta_1 =
0)>0$ to conclude that this Markov chain is {\it ergodic}: it
admits a unique stationary distribution, say $\pi$ (see,
e.g.,\cite{Tak}) and converges to it in the total variation norm,
for any initial value $\widetilde{V}_0\in {\Z}_+$. Further, note
that  $\pi$ is a Poisson distribution with parameter $c/p$ and
that the Markov chain is geometrically ergodic, i.e. its
distribution converges to the stationary one geometrically fast --
see Lemma 1 below.

Let a random variable $\widetilde{V}$ have the distribution $\pi$
and do not depend on anything else. Then random variable $B_1(\widetilde{V},p)$ has
a Poisson distribution with parameter $c$ and ${\mathbf P} (B_1(\widetilde{V},p)=1)=ce^{-c}.$

\begin{theorem}
\label{t_1} If, for some $c>0$,
\[
\mu(q) = \min\left\{\displaystyle\frac{c}{q},1\right\}
\]
then the stochastic system described by the 2-dimensional Markov chain
\begin{eqnarray*}
\label{lw} \left\{
      \begin{array}{rl}
q_{n+1} = q_n - {\bf I} (B_n(v_n,p)=1) + \xi_{n},\\
v_{n+1} =  v_n - B_n(v_n,p) + \widehat{B}_n(q_n -
v_n+\xi_n,\mu(q_n))
\end{array}
          \right.
\end{eqnarray*}
 is stable if $\lambda < ce^{-c}$
and unstable if $\lambda
> ce^{-c}$.
\end{theorem}

\begin{remark} The origin  $(0,0)$ is achievable from all other states $(x,y)\in \Z^2_+$ of the
2-dimensional Markov chain $(v_n,q_n)$. Since ${\mathbf P}
(\xi_1=0)>0$, the Markov chain is aperiodic. Therefore, if the
Markov chain is positive recurrent, it is also ergodic: there is a
stationary distribution, say $\phi$, on $\Z^2_+$ such that, for
any initial value $(q_0,v_0)$, the distribution of $(q_n,v_n)$
converges to $\phi$ in the total variation norm:
$$
\sup_{A\subset \Z^2_+} |{\mathbf P} ((q_n,v_n)\in A) - \phi (A)| \to 0, \quad n\to\infty.
$$
\end{remark}

\begin{remark}
The stability condition attains its maximum value $e^{-1}$ for
$c=1$ and doesn't depend on parameter $p.$
\end{remark}

\begin{remark} More or less straightforward modifications of the proof of Theorem 1
lead to the following results.
If, instead of \eqref{star}, we assume that
\[
q \cdot \mu(q) \to \infty, \text{ as } q\to\infty,
\]
then the system is similar to the classical ALOHA, which is
unstable. On the other hand, if we assume that
\[
q\cdot \mu(q) \to 0, \text{ as } q\to\infty,
\]
then the system is unstable again as it does not have enough energy
to transmit the ``old''  messages.
\end{remark}

\section{Auxiliary lemma}

We summarise a number of simple, but important properties in the following lemma.

\begin{lemma}\label{lem1} Let $\{Z_n\}$ be an i.i.d. sequence of non-negative integer-valued
random variables with finite mean ${\mathbf E} Z_0$. Assume it does not depend on i.i.d. random
variables $\{U_{n,i}, -\infty < n < \infty, i\ge 1\}$
having the uniform-$(0,1)$ distribution. \\
(I) For any initial value $W_0$ and for any $p\in (0,1)$, the sequence
\begin{equation}\label{WW}
W_{n+1} = W_n - B_n(W_n,p) + Z_n \equiv D_n(W_n,1-p) + Z_n
\end{equation}
is ergodic.\\
(II) For $m\le n$ define operators
$$
D_{m:n}(k,1-p) = D_n (D_{n-1}(\ldots (D_{m+1}(D_m(k,1-p),
1-p)\ldots , 1-p)
$$
and note that random variable $D_{m:n}(k,1-p)$ has a binomial distribution
with parameters $k$ and $(1-p)^{n-m+1}$. Then a stationary sequence $\{W^{(n)}\}$
given by
\begin{equation}\label{defW}
W^{(n)}=Z_{n-1} + \sum_{j=1}^{\infty} D_{(n-j-1):(n-2)}(Z_{n-j-1},1-p)
\end{equation}
has finite mean ${\mathbf E} W^{(n)}={\mathbf E} Z_0/p$ and forms
a stationary solution to recursive equation \eqref{WW},
$$
W^{(n+1)} = D_n(W^{(n)},1-p)+Z_n.
$$
(III) Let $Q$ be the distribution of $W^{(0)}$. Then, for any initial value $W_0$,
\begin{equation}\label{coupling}
{\mathbf P} (W_n=W^{(n)}) = {\mathbf P} (W_l=W^{(l)}, \ \ \mbox{for all} \ \ l\ge n) \to 1
\end{equation}
and, in particular,
$$
\sup_{A\subset \Z^2_+} |{\mathbf P} (W_n\in A)- Q(A)| \le
{\mathbf P} (W_n\ne W^{(n)}) \to 0, \ \
\mbox{as} \ \  n\to \infty.
$$
In particular, if $W_0$ and $Z_0$ have finite exponential moments, then
$$
{\mathbf P} (W_n\ne W^{(n)})\le K_1e^{-K_2n},
$$
for some $K_1,K_2>0$ and for all $N\ge 0$.\\
(IV) For any $n\ge 0$,
$$
W_n \le W_0 + W^{(n)} \quad \mbox{a.s.}
$$
and
$$
{\mathbf E} W_n \le {\mathbf E} W_0 + {\mathbf E} Z_0/p.
$$
(V) Let $\{\widetilde{Z}_n\}$ be any other sequence of non-negative and integer-valued random variables,
such that $0\le \widetilde{Z}_n\le Z_n$ a.s., for all $n$. Consider a recursion
$$
\widetilde{W}_{n+1} = D_n(\widetilde{W}_n,1-p)+\widetilde{Z}_n
$$
with integer initial value $0\le \widetilde{W}_0\le W_0$. Then
$$
\widetilde{W}_n \le W_n \quad \mbox{a.s., for all} \ \ n\ge 0.
$$
(VI) In particular, when $Z_0$ is a Poisson random variable with parameter $c$, formula \eqref{defW}
and the splitting and composition theorems imply that every $W^{(n)}$ has Poisson distribution with
parameter $c/p$.
\end{lemma}

The {\sc proof} of all statements of the lemma is straightforward.  First, we verify
that the right-hand side of \eqref{defW} is finite. This follows from finiteness
of ${\mathbf E} Z_0$ since expectation of the sum of non-negative random variables always equals to the
sum of expectations, and
\begin{eqnarray*}
{\mathbf E} Z_{n-1}+\sum_{j=1}^{\infty} {\mathbf E} (D_{(n-j-1):(n-1)}(Z_{n-j-1},1-p))
&=& \sum_{j=0}^{\infty} {\mathbf E} Z_0 \cdot (1-p)^j\\
&=&
{\mathbf E} Z_0 \cdot \frac{1}{p} <\infty.
\end{eqnarray*}
Therefore, the sum in the right-hand side of \eqref{defW} must contain only finite number of non-zero
elements. Further, if ${\mathbf E} e^{K_0Z_0}$ is finite for some $K_0>0$, then, for $K\in (0,K_0)$,
$$
{\mathbf E} e^{KW^{(0)}} = \prod_{n=0}^{\infty} {\mathbf E} \left(1+(e^K-1)(1-p)^n\right)^{Z_0}
$$
where all terms in the right product are finite and, as $K\downarrow 0$ and uniformly in $n$,
$$
{\mathbf E} \left(1+(e^K-1)(1-p)^n\right)^{Z_0} = 1+(1+o(1))K(1-p)^n{\mathbf E} Z_0
$$
so the product is finite too. Indeed, for any $h>0$,
$$
1+hZ_0 \le (1+h)^{Z_0}\le e^{hZ_0} \le 1+hZ_0+\frac{h^2Z_0^2}{2}e^{hZ_0}\le
1+hZ_0+\frac{h^2Z_0^4}{4}+\frac{h^2e^{2hZ_0}}{4}
$$
and then
$$
1+h{\mathbf E} Z_0\le {\mathbf E} (1+h)^{Z_0} \le 1+h{\mathbf E} Z_0 + O(h^2),\quad \mbox{as}\quad h\downarrow
0.
$$
Then one can take $h=(e^K-1)(1-p)^n$ and note that $e^K-1=K(1+o(1))$, as $K\downarrow 0$.

Next, at any time $n$, we have a set of $W_n$ elements, each of which disappears at time $(n+1)$
with
probability $p$, independently of everything else, while $Z_n$ new elements arrive. Then the original $W_0$
elements disappear in a (random) finite time, say $T$, and after that all $W_n$ include only new items that arrive
after time 0. Further, if ${\mathbf E} e^{K_0Z_0}$ is finite for some $K_0>0$, then, for
$k$ sufficienly small and uniformly in $n$,
\begin{eqnarray*}
{\mathbf P}(T>n) &=& {\mathbf P} (D_{0:(n-1)} (Z_0,1-p)\ge 1)\\
&=&
{\mathbf E} \left(1-(1-(1-p)^n)^{Z_0}\right)\\
&\le &
{\mathbf P} (Z_0>Kn) + 1-(1-(1-p)^n)^{Kn}\\
&\le &
\widetilde{K}e^{-Kn} + (1-p)^nKn (1+o(1))\\
&=& o\left(e^{-rn}\right),
\end{eqnarray*}
for some $\widetilde{K}$ and for any $r< \min \left( K, \log \frac{1}{1-p} \right).$
Similar observations hold for $W^{(n)}$: all elements of $W^{(0)}$ leave by finite time, say $\widehat{T}$, which has a finite exponential moment if $Z_0$ does.
Therefore $W_n$ and $W^{(n)}$ coincide after time $\max (T,\widehat{T})\le T+\widehat{T}$ which
also has a finite exponential moment.

\section{Proof of Theorem \ref{t_1}}

We start with the proof of stability. Let $X_n=\left(q_n,
v_n\right)$.

By the generalised Foster criterion (see, e.g.,
\cite{Tak}), it is enough  to find a test (``Lyapunov'') function
$L(x)\ge 0$ and a (sufficiently large) positive integer $N_0$ such
that the set
\[
D=\{ x\in {\Z}_+^2 \ : \ L(x)\leq N_0\},
\]
is compact
and, for
a properly chosen positive and bounded from above integer-valued function $g(x)$,
 the Markov chain has a bounded mean drift on $D$,
\begin{eqnarray}
\label{ddd}
\sup_{x\in D} \E\left(L\left(X_{g\left(X_0\right)}\right)-L\left(X_0\right)\ |\
X_0=x\right)\equiv K <\infty,
\end{eqnarray}
and the drift is uniformly negative on the complement $D^c$
of $D$: for any $x\in D^c$,
\begin{eqnarray}
\label{ccc}
\E\left(L\left(X_{g\left(X_0\right)}\right)-L\left(X_0\right)\ | \
X_0=x\right) \le -\varepsilon
\end{eqnarray}
where $\varepsilon$ is a fixed positive constant. We choose the test
(``Lyapunov'') function $L$ by
$$
L(x) = q+v, \quad \mbox{for} \quad x=\left(q,v\right)\in {\Z}_+^2,
\ v\leq q,
$$
and then a proper choice of function $g(x)$  leads to the
stability result.

Note that since
\begin{eqnarray*}
\E\left(q_1-q_0\ |\ (q_0, v_0)\right) = \lambda -
\P\left(B_0(v_0, p)=1 \ | \ (q_0, v_0)\right)\leq 1
\end{eqnarray*}
and
\[
\E\left(v_1-v_0 \ | \  (q_0, v_0)\right)=  - p v_0 +
(q_0-v_0+\lambda) \ {\bf I} (q_0\leq c)
\]
\[
+c\ \frac{q_0-v_0+\lambda}{q_0}\  {\bf I} (q_0> c)\leq 2
(q_0-v_0+\lambda)
\]
the drift $(\ref{ddd})$ is bounded from above on the set $D$, for any choice
of constant $N_0$.

Recall that $v_0\leq q_0$, so if we choose $N_0$ such
that $\frac{N_0}{2}>c$, then $q_0> c$ for $X_0\in D^c.$ Therefore,
the drift $(\ref{ccc})$ is uniformly negative for $g(x)=1$ for all
$(q_0,v_0)\in D^c,$ such that
$v_0>\frac{c+\lambda+1}{p}.$

Now we complete the stability proof, with
choosing sufficiently large values of $N_0$ and $k$ such that,
with $g(x)\equiv k$, inequality \eqref{ccc} uniformly holds for all $x=(q_0,v_0)\in D^c$ such that $v_0\le \frac{c+\lambda+1}{p}$.

Let $\{Y_n\}$ be i.i.d. random variables that do not depend on $\{\xi_n\}$ and have distribution
$$
{\mathbf P} (Y_n>x) = \sup_{q\ge c} {\mathbf P} (B_1(q,\mu
(q))>x).
$$
By Chernoff's inequality, for $x\ge c$ and any $\alpha >0$, the right-hand side of the latter equality
does not exceed
$$
\sup_{q\ge c}{\mathbf E} e^{\alpha B_1(q,c/q)}e^{-\alpha x}
=
\sup_{q\ge c}\left(1+(e^{\alpha}-1)\frac{c}{q}\right)^q e^{-\alpha x}\equiv C_1 e^{-\alpha x}
$$
where $C_1$ is finite since $\left(1+(e^{\alpha}-1)\frac{c}{q}\right)^q \to \exp (c(e^{\alpha}-1))
<\infty$ as $q\to\infty$. Therefore the distribution of $Y_1$ is proper and, moreover, has a finite
exponential moment.\\
Let now
\begin{equation}\label{ZYW}
Z_n = Y_n+\xi_n \quad \mbox{and} \quad W_{n+1}=D_n(W_n,1-p)+Z_n,
\quad \mbox{with} \quad W_0=v_0.
\end{equation}
Then, by property (V)
of Lemma 1,
$$
{\mathbf E} (v_n \ | \ v_0 ) \le {\mathbf E} (W_n \ | \ W_0) \le v_0+ C_2
$$
where
$$
C_2 = {\mathbf E} (Y_1+\xi_1)/p \equiv ({\mathbf E} Y_1 + \lambda )/p.
$$

Next, choose $\delta >0$ such that $\lambda + \delta < ce^{-c}$ and
consider the auxiliary Markov chain $\widetilde{V}_n$. Since it is ergodic, for
given $\delta>0$ and $R:= (c+\lambda +1)/p$ (we assume $R$ to be an integer), there
exists $l\in {\Z}_+$ such that, for any $n\geq l$,
\begin{equation}
\label{gb}
\sup_{\widetilde{V}_0\leq
R}\left|\P(B_n(\widetilde{V}_n,p)=1) -
ce^{-c}\right|<\frac{\delta}{3}.
\end{equation}
Then we choose $k>l$ such that
\begin{equation}\label{vareps}
-\varepsilon := k\lambda +C_2 - (k-l) (ce^{-c}-\delta) <0.
\end{equation}
Now choose $C_3$ so large that, for the sequence $\{W_n\}$ defined in \eqref{ZYW} and with
initial value $W_0=R$,
the probability of event
$$
A(C_3) = \{W_i \le C_3 \ \ \mbox{for all}  \ \ 0\le i\le k\}
$$
is not smaller than $1-\delta/3$.

The famous Poisson theorem allows the following simple extension:
as $q\to\infty$, the distribution of random variable
$B_1(q-v+\xi_1, c/q)$ converges in the total variation norm to the Poisson distribution 
with parameter $c$, uniformly in all $v\in \{0,1,\ldots,C_3\}$.
Therefore, one can choose $\widehat{q}$ such that the total
probability distance between the distribution of
$B_1(q-v+\xi_1,c/q)$ and Poisson distribution with parameter $c$
does not exceed $\delta/3$, for all $q\ge \widehat{q}$.

 Finally, let $N_0=\widehat{q}+k+R$. For any $(q_0,v_0)$
with $q_0+v_0\ge N_0$ and $v_0\le R$, we have $q_i\ge \widehat{q}$
for all $0\le i\le k$, since $q_{i+1}\ge q_i-1$ a.s. Then we may
put all bounds together to conclude that, under the assumptions
from above, for any $i=l+1,\ldots,k$,
$$
|{\mathbf P} (B_i(v_i,p)=1\ | \ (q_0,v_0))-ce^{-c}|\le \delta
$$
and, therefore, \eqref{ccc} holds. Indeed, with $x=(q_0,v_0)$ where $q_0+v_0\ge N_0$
and $v_0\le R$, we have
\begin{eqnarray*}
{\mathbf E} (L(X_k)-L(X_0) \ | \ X_0=x)  &=&
{\mathbf E} (q_k-q_0 \ | \ X_0=x) + {\mathbf E} (v_k-v_0 \ | \ X_0=x) \\
&=&
C_2 + k\lambda - \sum_{i=1}^k {\mathbf P} (B_i(v_i,p)=1 \ | \ X_0=x)\\
&\le &
 C_2 + k\lambda - \sum_{i=l+1}^k {\mathbf P} (B_i(v_i,p)=1 \ | \ X_0=x) \le
 -\varepsilon
 \end{eqnarray*}
 where $\varepsilon$ is from \eqref{vareps}.
This completes the proof of the
stability part of the theorem.

The proof of the second part of the theorem is based on Theorem 2.1 from \cite{F_D}.
Assume one may define a nonnegative function $\widetilde{L}(x)$  in
such a way that, for a properly chosen bounded positive
integer-valued function $\widetilde{g}(x)$ and for any initial
value $X_0=x$ with $\widetilde{L}(x)\geq \widetilde{N}_0,$  for
some fixed $\widetilde{N}_0>0,$  the Markov chain has a positive
drift
\begin{equation}
\label{c1} \E\left(\Delta_{\widetilde{g}(x)} {\bf
I}\left(\Delta_{\widetilde{g}(x)}\leq
M\right) \ | X_0=x\right)\geq\widetilde{\varepsilon}
\end{equation}
where
$\Delta_{\widetilde{g}(x)}=\widetilde{L}\left(X_{\widetilde{g}(x)}\right)-\widetilde{L}\left(x\right)$ given $X_0=x$, and the sequence
\begin{equation}
\label{c3} \left(\Delta_{\widetilde{g}(x)}^{-}\right)^2=\left(
\min\left\{0,\Delta_{\widetilde{g}(x)}\right\}\right)^2
\end{equation}
is uniformly integrable in all $x$ such that $\widetilde{L}(x)\geq
\widetilde{N}_0.$ Here $\widetilde{\varepsilon}$ and $M$ are fixed
positive constants. By Theorem 2.1 from \cite{F_D}, under
conditions $(\ref{c1})$-$(\ref{c3})$ imply that, for any $x\in
{\Z}_+^2$,
\[
P\left(\widetilde{L}(X_n)\to\infty \ | X_0=x \right)=1.
\]

We will apply this Theorem for the following choice of functions $\widetilde{L}$
and $\widetilde{g}$: for $x=(q,v)$,
we let $\widetilde{L}(x)=q$, and $\widetilde{g}(x)=1$ if $v > \widetilde{R}$
and $\widetilde{g}(x)=\widetilde{k}$ if $v\le \widetilde{R}$. We choose $\widetilde{R}$
and $\widetilde{k}$ below. Before that, we like to mention that our choice
of functions leads to redundancy of condition \eqref{c3} since $q_{i-1}\ge q_i-1$ a.s..
Also, since $\lambda ={\mathbf E} \xi_1$ is finite, the increments
$\Delta_{\widetilde{g}(x)}$ are uniformly integrable in $x$ and, therefore, condition
\eqref{c1} is equivalent to
\begin{equation}\label{c4}
{\mathbf E} \left( \Delta_{\widetilde{g}(x)} \right) \ge \widetilde{\varepsilon},
\end{equation}
for some positive $\widetilde{\varepsilon}$.

Let $\widetilde{R}$ be any positive integer such that
$$
\gamma := \sup_{n\ge \widetilde{R}}  {\mathbf P} (B_1(n,p)=1) < \lambda.
$$
Since ${\mathbf P} (B_1(n,p)=1) \to 0$ as $n\to\infty$, one can
always find such $\widetilde{R}$. Then, indeed, for
$v_0>\widetilde{R}$
$$
{\mathbf E} (q_1-q_0 \ | \ X_0=(q_0,v_0)) \ge \lambda - \gamma >0.
$$
Assume now that $\widetilde{R} >1$ and consider the case where
$q_0$ is large and $v_0 \le \widetilde{R}$. On the essence, we
repeat the scheme of the last part of the proof of stability, but
with considering the opposite inequalities.

First, we choose $0<\widetilde{\delta}< \lambda -ce^{-c}$ and find then $\widetilde{l}$
such that, for all $n\ge \widetilde{l}$,
\begin{equation}
\label{gb2}
\sup_{\widetilde{V}_0\leq
\widetilde{R}}\left|\P(B_n(\widetilde{V}_n,p)=1) -
ce^{-c}\right|<\frac{\widetilde{\delta}}{3}.
\end{equation}
Then we choose $\widetilde{k}>\widetilde{l}$ such that
$$
\widetilde{\varepsilon} := \widetilde{k}\lambda 
-\widetilde{l} -
(\widetilde{k}-\widetilde{l}) (ce^{-c}+\widetilde{\delta}) >0.
$$
Then we choose $\widetilde{C}$ so large that, for the sequence $\{W_n\}$ defined in \eqref{ZYW}
and with initial value $W_0=\widetilde{R}$,
the probability of event
$$
A(\widetilde{C}) = \{W_i \le \widetilde{C} \ \ \mbox{for all}  \ \ 0\le i\le n\}
$$
is not smaller than $1-\widetilde{\delta}/3$.\\
Finally we choose $\widehat{q}$ such that the total probability
distance between the distribution of $B_1(q-v+\xi_1,c/q)$ and
Poisson distribution with parameter $c$ does not exceed
$\widetilde{\delta}/3$, for all $q\ge \widehat{q}$ and $v\le
\widetilde{C}$.

Then, putting altogether, we obtain that if $q_0>\widehat{q}+\widetilde{k}$, then
$$
{\mathbf E} (q_{\widetilde{k}}-q_0 \ | \ X_0= (q_0,v_0))\ge \widetilde{\varepsilon}
$$
if $v_o\le \widetilde{R}$. This completes the proof of the second part of the theorem.

\section{Conclusion}
We studied stability of the random access system under stochastic
energy harvesting for a model with infinite number of transmitting nodes. We assumed
that the harvesting intensity of individual power
batteries depends on the number of active users. General stability
and instability conditions  have been obtained.

The model introduced in this paper assumes an individual
energy supply mechanism for each users. One may comment that, in
existing wireless systems, the energy harvesting mechanism is such
that each user has an individual storage element (battery or
capacitor) which is charged by the common energy mechanism (see,
for instance,
 \cite{XRC}).
 If the number
of users is relatively low, then the amount of energy arriving at
a single user does not depend on the total number of users that
are storing the energy. However, if the number of users becomes
sufficiently large, the energy arriving at an individual user
starts to decrease. This effect was discussed in \cite{Can},
\cite{Tong}. Thus, the models introduced in this paper reflects
the features of wireless systems with an energy harvesting
mechanism that has a single base station and a large number of
user devices.

\bigskip

\begin{otherlanguage}{english}

\end{otherlanguage}


\begin{thebibliography}{1}

\bibitem{XRC}  Xun Zhou, Rui Zhang, Chin Keong
Ho, {\it Wireless Information and Power Transfer in Multiuser OFDM
Systems},  IEEE Transactions on Wireless Communications, {\bf
13}:4 (2014), 2282-2294.



\bibitem{Ephr}
J. ~Jeon and A. ~Ephremides, {\it The stability region of random
multiple access under stochastic energy harvesting}, Proceedings
of the IEEE International Symposium on Information Theory (ISIT),
(2011), 1796-1800.

\bibitem{sib}
S.~Foss, D.~Kim, A.~Turlikov,  {\it Random multiple access system
with a common energy supply}, XIV International symposium on
problems of redunduncy in information and control system, (2014),
 39-42.

\bibitem{Abr}
N.~Abramson, {\it Development of the ALOHANET}, IEEE Trans. Info.
Theory, {\bf 31} (1985),  119-123.

\bibitem{FKT}
S.~Foss, D.~Kim, A.~Turlikov,  {\it On the Models of Random
Multiple Access with Stochastic Energy Harvesting}, The 7th
International Congress on Ultra Modern Telecommunications and
Control Systems, (2015).


\bibitem{Tak}
S.~Foss and T.~Konstantopoulos, {\it An overview of some
stochastic stability methods}, Journal of Operation Research
Society Japan, {\bf 47}:4 (2004), 275-303.

\bibitem{Can}
B.~L. Cannon, J.~F. Hoburg,  D.~D. Stancil, S.~C. Goldstein, {\it
Magnetic Resonant Coupling As a Potential Means for Wireless Power
Transfer to Multiple Small Receivers}, IEEE Transactions on Power
Electronics, {\bf 24}:7 (2009), 1819-1825.

\bibitem{Tong}
B. Tong , Z. Li , G. Wang and W. Zhang, {\it How Wireless Power
Charging Technology Affects Sensor Network Deployment and
Routing}, Proc. IEEE 30th Int'l Conf. Distributed Computing
Systems (ICDCS), (2010),  438-447.


\bibitem{F_D}
S.~Foss and D.~Denisov, {\it On transience conditions for Markov
chains}, Siberian Mathematical Journal, {\bf 42}:2 (2001),
364-371.

\end{thebibliography}
\end{document}